\input epsf
\documentclass[12pt]{article}
\usepackage{epsfig}
\usepackage{amssymb}
\usepackage{graphicx}
\newcommand{\gil}[1]{\lfloor #1\rfloor}
\newcommand{\lig}[1]{\lceil #1\rceil}

\newcommand{\be}{\begin{equation}}
\newcommand{\ee}{\end{equation}}
\newcommand{\qed}{\hfill $\square$\vskip .2cm}

\newtheorem{remark}{Remark}[section]
\newtheorem{example}{Example}[section]

\newcommand{\p}{{\mathbb P}}
\newcommand{\G}{{\mathbb G}}
\newcommand{\n}{{\mathbb N}}
\newcommand{\C}{{\mathbb C}}

\newcommand{\E}{{\mathbb E}}

\renewcommand{\S}{{\mathbb S}}

\newcommand{\B}{{\cal B}}

\newcommand{\R}{{\mathbb R}}

\def\<{\langle}
\def\>{\rangle}

\newtheorem{prop}{Proposition}[section]

\newtheorem{lem}{Lemma}[section]
\newtheorem{theorem}{Theorem}[section]

\newcommand{\sect}[1]{\section{#1}\setcounter{equation}{0}}

\begin{document}
\title{Occupation laws for some time-nonhomogeneous Markov chains}

\author{Zach Dietz$^1$\ \ and \ \ Sunder Sethuraman$^2$}

\thispagestyle{empty}

 \maketitle
 \abstract{
We consider finite-state time-nonhomogeneous Markov chains whose
 transition matrix at time $n$ is $I+G/n^\zeta$ where $G$ is a ``generator''
 matrix, that is $G(i,j)>0$ for $i,j$ distinct, and $G(i,i)= -\sum_{k\neq
 i}G(i,k)$, 
and $\zeta>0$ is a strength parameter.
In these chains, as time grows, the positions are less and
 less likely to change, and so form simple models of age-dependent
time-reinforcing schemes.  These chains, however, exhibit some 
different, perhaps unexpected, occupation behaviors depending on parameters.

 Although it is shown, on the one hand, that the position at time $n$ converges to
 a point-mixture for all $\zeta>0$, on the other hand, the average
 occupation vector up to time $n$, when
 variously $0<\zeta<1$, $\zeta>1$
 or $\zeta=1$, is seen to converge to a constant, a point-mixture, or a
 distribution $\mu_G$ with no atoms and
 full support on a simplex respectively, as $n\uparrow \infty$.  This
 last type of limit
can be interpreted as
 a sort of ``spreading'' between the cases $0<\zeta<1$ and
 $\zeta>1$.

In particular, when $G$ is
 appropriately chosen, intriguingly, $\mu_G$ is a
 Dirichlet distribution, reminiscent of results in P\'olya urns.
}

 \vskip .1cm
 \thanks{Research supported in part by NSA-H982300510041 and NSF-DMS-0504193 
 \\
 {\sl Abbreviated title}: Occupation laws for nonhomogeneous Markov chains.  \\
 {\sl AMS (2000) subject classifications}: Primary 60J10; secondary
 60F10.}\\
$^1$ Department of Mathematics, Tulane
 University, 6823 St. Charles Ave., New Orleans, LA \ 70118;
 zdietz@math.tulane.edu. \\
$^2$ Department of Mathematics, Iowa State University,
 396 Carver Hall, Ames, IA  \ 50011; sethuram@iastate.edu.

\sect{Introduction and Results}

In this article, we study laws of large numbers (LLN) for
a class of finite space
time-nonhomogeneous Markov chains
where, as time increases, positions are less likely to change.  Although these chains feature simple age-dependent
time-reinforcing dynamics, some different, perhaps unexpected, 
LLN occupation behaviors
emerge depending on parameters.  A specific case, as in Example 1.1, was first introduced in Gantert
\cite{Gantert} in connection with analysis of certain simulated
annealing LLN phenomena. 

\begin{example}\rm Suppose there are only two states $1$ and $2$, and
  that
the chain moves between the two locations in the following way:  At large
times $n$, the chain switches places with probability $c/n$, and stays
put with complementary probability 
$1-c/n$ for $c>0$.  The chain, as it ages, is less inclined to leave its spot,
but nonetheless switches infinitely often.  One can see the
probability of being in state $1$ tends to $1/2$ regardless of the
initial distribution.  One may ask, however, how the average location,
or frequency, of state $1$ behaves asymptotically.  For this example,
it was shown in \cite{Gantert} and Ex. 7.1.1. \cite{Winkler}, 
perhaps surprisingly, that any LLN limit could not be
a constant, or even converge in probability, without further identification.  However, a quick consequence of our results is that the average
occupation limit of state $1$ converges weakly to the Beta$(c,c)$
distribution
(Theorem \ref{Dirichletthm}).  
\end{example}

More specifically, we 
consider a general version of this scheme with $m\geq 2$ possible
locations, and moving and staying probabilities $G(i,j)/n^\zeta$ and
$1-\sum_{k\neq i}G(i,k)/n^\zeta$ from $i\rightarrow j\neq i$ and
$i\rightarrow i$ respectively at time $n$ where $G=\{G(i,j)\}$ is an
$m\times m$ matrix and
$\zeta>0$ is a strength parameter.  
After observing the location probabilities tend to a
distribution which depends on $G$, $\zeta$, and initial probability
$\pi$ when $\zeta>1$, but does not depend on $\zeta$ and $\pi$ when
$\zeta\leq 1$ (Theorem
\ref{marginalprop}), the results on the average occupation vector limit separate roughly 
into three cases
depending on whether
 $0<\zeta<1$, $\zeta=1$, or $\zeta>1$.

When $0<\zeta<1$, following \cite{Gantert}, 
the average occupation is seen to converge to a
constant in probability; and when more specifically
$0<\zeta<1/2$, this convergence is proved to be a.s.
When $\zeta>1$, as there are only a finite number of
switches, 
the position eventually stabilizes and the average
occupation converges to a mixture of point masses (Theorem \ref{theorem1}).

Our main results are when $\zeta=1$.  In this case, 
we show the
average occupation converges to a non-atomic distribution $\mu_G$, with full
support on a simplex, identified by its moments (Theorems
\ref{existencethm} and \ref{structurethm}).  When, in particular, $G$
takes form $G(i,j) = \theta_j$ for all $i\neq j$, that is when the 
transititions into a state $j$ are constant, $\mu_G$ takes the form of 
a Dirichlet distribution with parameters
$\{\theta_j\}$ (Theorem \ref{Dirichletthm}).
The proofs of these statements
follow by the method of moments, and some 
surgeries of the paths.

The heuristic is that when $0<\zeta<1$ the chance of switching
is strong and sufficient mixing leads to constant limits, 
but when
$\zeta>1$ there is little movement giving point-mixture limits. The
case $\zeta =1$ is the intermediate ``spreading'' situation leading to
non-atomic limits.  For example, with respect to Ex. 1.1, when the
switching probability at time $n$ is $c/n^\zeta$, the Beta$(c,c)$ limit when $\zeta=1$
interpolates, as $c$ varies on $(0,\infty)$, between the point-mass at
$1/2$, the frequency limit of state $1$ when
$0<\zeta<1$, and the fair mixture of point-masses at $0$ and $1$,
the limit when $\zeta>1$ and starting at random (cf. Fig. 1).


\begin{figure}
\begin{center}
\psfig{figure=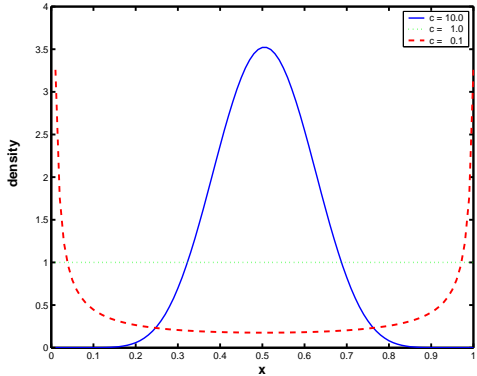,height=2in,width=2.5in}

\caption{Beta$(c,c)$ occupation law of state $1$ in Ex. 1.1.}
\end{center}
\end{figure}

In the literature, there are only a few results on LLN's for
time-nonhomogeneous Markov chains, often related to simulated
annealing and Metropolis
algorithms which can be viewed in terms of a generalized model
where $\zeta=\zeta(i,j)$ is a non-negative function. 
These results relate to the
case ``$\max\zeta(i,j)<1$'' when the LLN limit is a constant \cite{Gantert}, Ch. 7 \cite{Winkler}, \cite{Gidas}.  See also Ch. 1 \cite{IT}, \cite{Liu},\cite{Liu-Yang}; and
texts \cite{Bremaud},
\cite{Isaacson},\cite{Iosifescu} for more on nonhomogeneous
Markov chains.  In this light, the non-degenerate limits
$\mu_G$ found here seem to be novel objects. In terms of simulated
annealing, these limits suggest a more complicated LLN picture at
the ``critical'' cooling schedule when $\zeta(i,j)=1$ for some pairs
$i,j$ in the state space.
 


%
%
%
The advent of Dirichlet limits, when $G$ is chosen appropriately,
seems of particular interest, given
similar results for limit color-frequencies in P\'olya urns
\cite{Athreya}, \cite{Gouet}, as it hints at an even larger role
for Dirichlet
measures in related but different ``reinforcement''-type models (see
\cite{KB}, \cite{Pemantle}, \cite{dMM}, and
references therein, for more on urn and reinforcement schemes).  
In this context, the set of
``spreading'' limits $\mu_G$ in Theorem \ref{existencethm}, in which Dirichlet
measures are but a subset, appears intriguing as well (cf. Remarks 1.4, 1.5 and Fig. 2).

In another vein, although different, Ex. 1.1 seems
not so far from the case of independent Bernoulli trials with success
probability $1/n$ at the $n$th trial.  For such trials much is
known about the spacings between successes, and connections to GEM
random allocation models and Poisson-Dirichlet measures \cite{Vervaat}, \cite{Arratia}, \cite{ABT},
\cite{ABT-book}, \cite{Pitman}, \cite{Pitman-book}. 

We also mention, in a different, neighbor setting, some interesting
but distinct
LLN's have
been shown for
arrays of time-{\it homogeneous}
Markov sequences where the transition matrix $P_n$ for the $n$th row
converges to a limit matrix $P$ \cite{Dob}, \cite{Hanen}, Section 5.3 \cite{Iosifescu}; see also
\cite{Miclo} which comments on some ``metastability'' concerns.

We now develop some notation to state results.  Let $\Sigma =
\{1,2,\ldots,m\}$ be a finite set of $m\geq 2$ points. We say a
matrix $M=\{M(i,j): 1\leq i,j\leq m\}$ on $\Sigma$ is a {\it generator}
matrix if $M(i,j)\geq 0$ for all distinct $1\leq i,j\leq m$, and
$M(i,i) = -\sum_{j\neq i}M(i,j)$ for $1\leq i\leq m$.  In particular, $M$ is
a generator with {\it nonzero entries} if $M(i,j)>0$ for $1\leq i,j\leq m$
distinct, and $M(i,i)<0$ for $1\leq i\leq m$.

To avoid technicalities, e.g. with reducibility, 
we work with the following matrices, \begin{eqnarray*} \G &=&\bigg\{G\in\R^{m\times m}: G\ \mbox{is
a generator matrix with nonzero entries}\bigg\}, \end{eqnarray*}
although extensions should be possible for a larger class.
 For
$G\in \G$, let $n(G,\zeta) = \lig{\max_{1\leq i\leq
m}|G(i,i)|^{1/\zeta}}$, and define for $\zeta>0$
$$P^{G,\zeta}_n \ =\ \left\{\begin{array}{rl} I & \ {\rm for  \ } 1\leq n \leq n(G,\zeta)\\
 I+G/n^\zeta & \ {\rm for \ } n\geq n(G,\zeta)+1\end{array}\right.$$
 where $I$ is the $m\times m$ identity
matrix.  Then, for all $n\geq 1$, $P^{G,\zeta}_n$ is ensured to be a
stochastic matrix.

Let $\pi$ be a distribution on $\Sigma$, and let $\p_\pi^{G,\zeta}$
be the (nonhomogeneous) Markov measure on the sequence space
$\Sigma^\n$ with Borel sets $\B(\Sigma^\n)$ corresponding to initial
distribution $\pi$ and transition kernels $\{P^{G,\zeta}_n\}$. That
is, with respect to the coordinate process, ${\bf
X}=\<X_0,X_1,\ldots\>$, we have $\p^{G,\zeta}_\pi(X_0 = i)=\pi(i)$
and the Markov property
$$\p^{G,\zeta}_\pi(X_{n+1} = j|X_0,X_1,\ldots,X_{n-1},X_n=i)
= P^{G,\zeta}_{n+1}(i,j)$$ for all $i,j\in \Sigma$ and $n\geq 0$. Our
convention then is that $P^{G,\zeta}_{n+1}$ controls
``transitions'' between times $n$ and $n+1$. Let also
$\E_\pi^{G,\zeta}$ be expectation with respect to
$\p^{G,\zeta}_\pi$.  More generally, $E_\mu$ denotes expectation
with respect to measure $\mu$.

Define the occupation statistic ${\bf Z}_n=\<Z_{1,n},\cdots,Z_{m,n}\>$ for $n\geq 1$
where
$$Z_{i,n}\ =\ \frac{1}{n}\sum_{k=1}^n 1_i(X_k)$$
for $1\leq i\leq m$. Then, ${\bf Z}_n$ is an element of the
$m-1$-dimensional simplex, $$\Delta_m \ = \ \bigg\{{\bf x}: \sum_{i=1}^m x_i
= 1, 0\leq x_i\leq 1 {\rm \ for \ } 1\leq i\leq m\bigg\}.$$


The first result is on convergence of the position of the process.
For $G\in \G$, let $\nu_G$ be the stationary distribution
corresponding to $G$ (of the associated continuous time
homogeneous Markov chain), that is the unique left eigenvector, with
positive entries,
normalized to unit sum, of the eigenvalue $0$.
\begin{theorem}
\label{marginalprop} For $G\in\G$, $\zeta>0$, and initial distribution $\pi$,
under $\p_\pi^{G,\zeta}$, \begin{eqnarray*} X_n \
\stackrel{d}{\longrightarrow}\ \nu_{G,\pi,\zeta} \end{eqnarray*}
where $\nu_{G,\pi,\zeta}$ is a probability vector on $\Sigma$
depending in general on $\zeta,\ G$, and $\pi$. When $0< \zeta\le
1$, $\nu_{G,\pi,\zeta}$ does not depend on $\pi$ and $\zeta$
and reduces to $\nu_{G,\pi,\zeta}=\nu_G$.
\end{theorem}

\begin{remark} \rm
For $\zeta>1$, with only finitely many moves, the convergence is a.s., and $\nu_{G,\pi,\zeta}$ is explicit when
$G=V_GD_GV_G^{-1}$ is diagonalizable with $D_G$ diagonal and
$D_G(i,i) = \lambda^G_i$, the $i$th eigenvalue of $G$, for $1\leq
i\leq m$. By calculation, $\nu_{G,\pi,\zeta} = \pi^t\prod_{n\geq
  1}P_n^{G,\zeta} = \pi^tV_G
D' V_G^{-1}$ with $D'$ diagonal and $D'(i,i)
=\prod_{n\geq n_0(G,\zeta)+1}(1+\lambda^G_i/n^\zeta)$.
\end{remark}

We now consider the cases $\zeta \neq 1$ with respect to average
occupation limits.  Let ${\bf i}$ be the basis vector
${\bf i}=\<0,\ldots,0,1,0,\ldots,0\>\in \Delta_m$ with a $1$ in the $i$th
component and $\delta_{\bf i}$ be the point mass at ${\bf i}$ for
$1\leq i\leq m$.

\begin{theorem}
\label{theorem1} Let $G\in\G$, and $\pi$ be an initial distribution.
Under
$\p^{G,\zeta}_\pi$, we have that
\begin{eqnarray*} {\bf Z}_n
 \ \longrightarrow \ \nu_G\  \end{eqnarray*}
converges in probability when $0<\zeta<1$; when more
specifically $0<\zeta<1/2$, this convergence is $\p^{G,\zeta}_\pi$-a.s.

However, when $\zeta>1$, under $\p^{G,\zeta}_\pi$, \begin{eqnarray*} {\bf Z}_n \
\stackrel{d}{\longrightarrow} \ \sum_{i=1}^m\nu_{G,\pi,\zeta}(i)
\delta_{\bf i} \ .\end{eqnarray*}
\end{theorem}

\begin{remark}\rm
Simulations suggest that actually a.s. convergence might hold also on the range
$1/2\leq \zeta<1$ (with worse convergence rates as
$\zeta\uparrow 1$).
\end{remark}

Let now $\gamma_1,\ldots,\gamma_m\geq 0$, be integers such that
$\bar{\gamma}=\sum_{i=1}^m \gamma_i\geq 1$. Define the list $A = \{a_i: 1\leq
i\leq \bar{\gamma}\}=
\{\underbrace{1,\ldots,1}_{\gamma_1},\underbrace{2,\dots,2}_{\gamma_2},
\dots, \underbrace{m,\dots,m}_{\gamma_m}\}$. Let
$\S(\gamma_1,\ldots,\gamma_m)$ be the $\bar{\gamma}!$ permutations
of $A$, although there are only
${\bar{\gamma}}\choose{\gamma_1,\gamma_2,\cdots,\gamma_m}$ distinct
permutations; that is, each permutation appears $\prod_{k=1}^m
\gamma_k !$ times.

Note also, for $G\in
\G$, being a generator matrix, all eigenvalues of $G$ have
non-positive real parts (indeed, $I+G/k$ is a stochastic matrix for
$k$ large; then, by Perron-Frobenius, the real parts of its
eigenvalues satisfy $-1\leq 1+{\rm Re}(\lambda^G_i)/k\leq 1$, yielding the
non-positivity), and
so the resolvent $(xI-G)^{-1}$ is well defined for $x\ge 1$.

\begin{theorem}
\label{existencethm} For $\zeta=1$, $G\in \G$, and initial
distribution $\pi$, we have under $\p^{G,\zeta}_\pi$ that $${\bf Z}_n\
\stackrel{d}
{\longrightarrow}\ \mu_G$$ where $\mu_G$ is a measure on the simplex
$\Delta_m$ characterized by its moments:  For $1\le i\le m$,
$$E_{\mu_G}\big(x_i\big) \ = \   \lim_{n\rightarrow \infty}
\E_{\pi}^{G,\zeta}\big(Z_{i,n}\big)\ = \ \nu_G(i), $$
and for integers $\gamma_1,\ldots,\gamma_m\geq 0$
when $\bar{\gamma}\geq 2$,
\begin{eqnarray*}E_{\mu_G}\bigg(x_1^{\gamma_1}\cdots
  x_m^{\gamma_m}\bigg)
&=&\lim_{n\rightarrow \infty}
\E_{\pi}^{G,\zeta}\bigg(Z_{1,n}^{\gamma_1}\cdots
Z_{m,n}^{\gamma_m}\bigg)\\
& =& \frac{1}{\bar{\gamma}} \sum_{\sigma\in
\S(\gamma_1,\ldots,\gamma_m)} \nu_{G}(\sigma_1)
\prod_{i=1}^{\bar{\gamma}-1}
\bigg(iI-G\bigg)^{-1}(\sigma_{i},\sigma_{i+1}). \end{eqnarray*}

\end{theorem}

 \begin{remark}\rm
 However, as in Ex. 1.1 and \cite{Gantert}, when $\zeta =1$ as above, ${\bf
   Z}_n$ cannot converge in probability (as the tail field $\cap_n \sigma\{X_n,X_{n+1},\ldots\}$ is trivial by
   Theorem 1.2.13 and Proposition 1.2.4 \cite{IT} and (\ref{tailtrivial}),
   but the limit distribution
   $\mu_G$ is not a point-mass by say Theorem \ref{structurethm} below).
This is in contrast to P\'olya urns where the color frequencies
   converge a.s.
 \end{remark}

We now consider a particular matrix under which $\mu_G$ is a
Dirichlet distribution. For $\theta_1,\ldots,\theta_m>0$,
define
\begin{eqnarray*} \Theta&=&\left[ \begin{array}{ccccc}
\theta_1-\bar{\theta}&\theta_2&\theta_3&\cdots&\theta_m\\
\theta_1&\theta_2-\bar{\theta}&\theta_3&\cdots&\theta_m\\
\vdots   &  \vdots   &\ddots&\cdots&\vdots\\
\theta_1 &\theta_2 &\theta_3&\cdots &\theta_m-\bar{\theta} \end{array}
\right] \end{eqnarray*} where $\bar{\theta}=\sum_{l=1}^m \theta_l$. It is
clear $\Theta\in \G$. Recall identification of the Dirichlet
distribution by its density and moments; see \cite{KBJ}, \cite{JS} for more on
these distributions. Namely, the Dirichlet distribution on the
simplex $\Delta_m$ with parameters $\theta_1,\ldots,\theta_m$
(abbreviated as Dir$(\theta_1,\ldots,\theta_m)$) has density
$$
\frac{\Gamma(\bar{\theta})}{\Gamma(\theta_1)\cdots\Gamma(\theta_m)}
\ x_1^{\theta_1-1}\cdots x_m^{\theta_m-1}. $$
   The moments
 with respect to integers $\gamma_1,\ldots,\gamma_m\geq 0$ with
 $\bar{\gamma}\geq 1$ are \begin{equation}
\label{dirichletmoments}
E\bigg(x_1^{\gamma_1}\cdots x_m^{\gamma_m}\bigg)\ =\
\frac{\prod_{i=1}^{m}\theta_i(\theta_i+1)\cdots(\theta_i+\gamma_i-1)}
{\prod_{i=0}^{\bar{\gamma}-1}(\bar{\theta}+i)}, \end{equation} where we take
$\theta_i(\theta_i+1)\cdots(\theta_i+\gamma_i-1)=1$ when
$\gamma_i=0$.

\begin{theorem}
\label{Dirichletthm} We have $\mu_\Theta  =   {\rm
Dir}(\theta_1,\ldots,\theta_m)$.
\end{theorem}

\begin{remark}\rm
Moreover, by comparing the first few moments in Theorem \ref{existencethm} with (\ref{dirichletmoments}), one can check $\mu_G$ is
not a Dirichlet measure for many $G$'s with $m\geq 3$.  However, when
$m=2$, then any $G$ takes the form of $\Theta$ with $\theta_1= G(2,1)$
and $\theta_2 = G(1,2)$, and so $\mu_G = {\rm Dir}(G(2,1),G(1,2))$.
\end{remark}

We now characterize the measures $\{\mu_G: G\in \G\}$ as
``spreading'' measures different from the limits
when $0<\zeta<1$ and $\zeta>1$.

\begin{theorem}
\label{structurethm} Let $G\in \G$. Then, (1) $\mu_G(U)>0$ for any
non-empty open set $U\subset \Delta_m$. Also, (2) $\mu_G$ has no atoms.
\end{theorem}



  \begin{figure}
  \begin{center}
  \psfig{figure=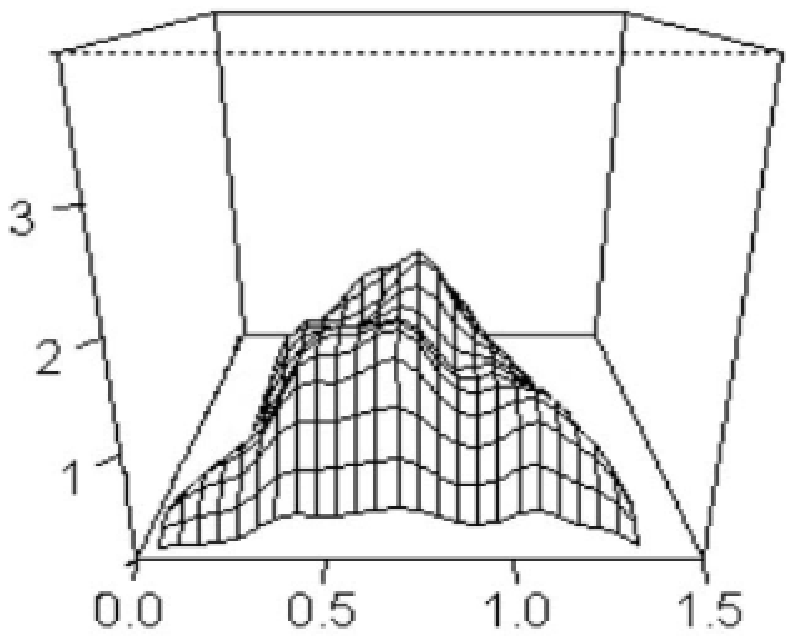,height=2in,width=2.5in}
  \psfig{figure=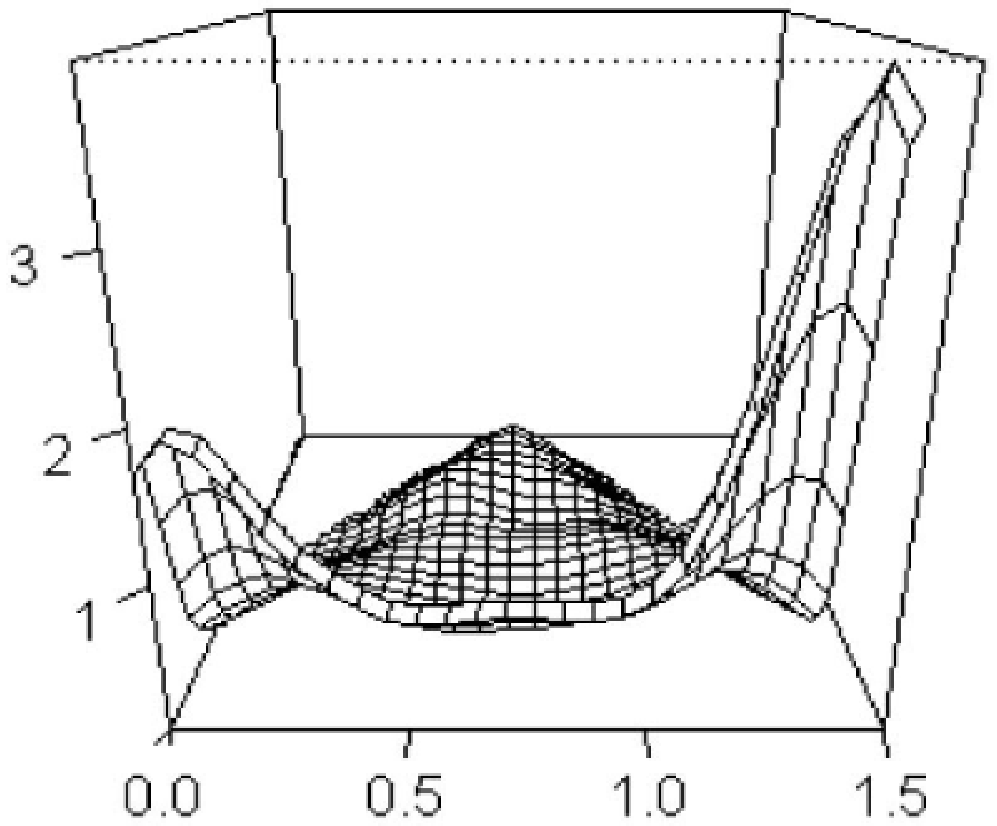,height=2.15in,width=2.5in}
  \caption{Empirical $\mu_G$ densities under $G_{\rm left}$ and
  $G_{\rm right}$ respectively.}
  \end{center}
  \end{figure}
\begin{remark}
\rm We suspect better estimates in the proof of Theorem
\ref{structurethm} will show $\mu_G$ is in fact mutually absolutely
continuous with respect to Lebesgue measure on $\Delta_m$.  Of
course, in this case, it would be of interest to find the density of
$\mu_G$. 
Meanwhile, we give two
histograms, found by calculating $1000$ averages, each on a run of
time-length $10000$ starting at random on $\Sigma$ at time $n(G,1)$
$(= 3,1 \ {\rm respectively})$,
 in Figure
2 of the empirical density when $m=3$ and $G$ takes forms
$$G_{\rm left} = \left[\begin{array}{rrr}-3&1&2\\
2&-3&1\\
1&2&-3\end{array}\right], \ \ {\rm and \ \ } G_{\rm right} =
\left[\begin{array}{rrr}-.4&.2&.2\\
.3&-.6&.3\\
.5&.5&-1\end{array}\right].$$ 
To help visualize plots, $\Delta_3$
is mapped to the plane by linear transformation
$f({\bf x}) = x_1f(\<1,0,0\>) + x_2f(\<0,1,0\>) +
x_3f(\<0,0,1\>)$ where $f(\<1,0,0\>)= \<\sqrt{2},0\>$,
$f(\<0,1,0\>)= \<0,0\>$ and $f(0,0,1)=
\sqrt{2}\<1/2,\sqrt{3}/2\>$. The map maintains a distance
$\sqrt{2}$ between the transformed vertices.
\end{remark}
%

%

We now comment on the plan of the paper.  The proofs of Theorems
\ref{marginalprop} and \ref{theorem1}, \ref{existencethm},
\ref{Dirichletthm}, and \ref{structurethm} (1) and (2) are in
sections 2,3,4, 5, and 6 respectively.  These sections do not depend
structurally on each other.

%
%
%

\sect{Proofs of Theorems \ref{marginalprop} and \ref{theorem1}}

We first recall some results for
nonhomogeneous Markov chains in the literature.  For a stochastic
matrix $P$ on $\Sigma$, define the ``contraction coefficient''
\begin{eqnarray}c(P) & =&
  \max_{x,y}\frac{1}{2}\sum_{z}\bigg|P(x,z)-P(y,z)\bigg|\nonumber\\
&=& 1- \min_{x,y}\sum_z \min\bigg\{P(x,z), P(y,z)\bigg\}
\label{contractioncoefficient}
\end{eqnarray}
The following is, for instance, Theorem 4.5.1 \cite{Winkler}.
\begin{prop}
\label{dobrushinprop} Let $X_n$ be a time-nonhomogeneous Markov
chain on $\Sigma$ connected by transition matrices $\{P_n\}$ with
corresponding stationary distributions $\{\nu_n\}$. Suppose
\begin{equation}
\label{dobeq}\prod_{n=1}^\infty c(P_n) = 0 \ \ {\it  and\ \ }
\sum_{n=1}^\infty \|\nu_n - \nu_{n+1}\|_{\rm Var}
<\infty.\end{equation} Then, $\nu = \lim_{n\rightarrow \infty}\nu_n$
exists, and, starting from any initial distribution $\pi$, we have
for each $k\in \Sigma$ that
$$\lim_{n\rightarrow \infty}P(X_n = k) \ =\  \nu(k).$$
\end{prop}

The following is stated in Section 2 \cite{Gantert} as a consequence of
results (1.2.22) and Theorem 1.2.23 in \cite{IT}.
\begin{prop}
\label{gantertprop} Given the setting of Proposition
\ref{dobrushinprop}, suppose (\ref{dobeq}) is satisfied, and 
$c_n = \max_{n_0\leq i\leq n}c(P_i)  <  1$ for all $n\geq n_0$ for
some $n_0\geq 1$.
Let $\pi$ and $f$ be any initial distribution, and function
$f:\Sigma\rightarrow \R$.  Then, we have convergence
 $$\frac{1}{n}\sum_{i=1}^n f(X_i) \
\rightarrow \ E_\nu[f]$$ 
in the following senses:

(i) In probability, when $\lim_{n\rightarrow \infty} n(1-c_n) =
\infty$.

(ii) a.s. when $\sum_{n\geq n_0} 2^{-n}(1-c_{2^n})^{-2}  <  \infty$.
\end{prop}

\vskip .1cm {\it Proof of Theorem \ref{marginalprop}.} We first
consider when $\zeta>1$. In this case there are only a finite number
of movements by Borel-Cantelli since $\sum_{n\geq 1}
\p^{G,\zeta}_\pi(X_n\neq X_{n+1}) \leq C\sum_{n\geq 1} n^{-\zeta}
<\infty$.  Hence there is a time of last movement $N<\infty$ a.s.
Then, $\lim X_n = X_N$ a.s., and, for $k\in \Sigma$, the limit distribution $\nu_{G,\pi,\zeta}$
is defined and given by
$\p_\pi^{G,\zeta}(X_N = k) =  \nu_{G,\pi,\zeta}(k)$.

When $0<\zeta\leq 1$, as $G\in \G$, by calculation with
(\ref{contractioncoefficient}),
$c(P^{G,\zeta}_n)= 1- C_{G}/n^{\zeta}$ for all $n\geq n_0(G,\zeta)$
large enough and 
a constant $C_G>0$.  Then,
\begin{equation}
\label{tailtrivial}\prod_{n\geq
1}c(P^{G,\zeta}_n)\ =\ \prod_{n\geq n_0(G,\zeta)}\bigg(1-\frac{C_G}{n^\zeta}\bigg)\ =\ 0.\end{equation}
Since for $n> n(G,\zeta)$,
$\nu_G^tP^{G,\zeta}_n=\nu_G^t(I-G/n^{\zeta}) = \nu_G^t$, the second
condition of Proposition \ref{dobrushinprop} is trivially satisfied,
and hence the result follows. \qed

\vskip .1cm

{\it Proof of Theorem \ref{theorem1}.}  When $\zeta>1$, as mentioned
in the proof of Theorem \ref{marginalprop}, there are only a finite
number of moves a.s., and so a.s. $\lim {\bf Z}_n = \sum_{k=1}^m 1_{[X_N =  k]}{\bf k}$ concentrates on basis
vectors $\{{\bf k}\}$.  Hence, as defined in proof of Theorem \ref{marginalprop}, 
$\p^{G,\zeta}_\pi(X_N = k)  = \nu_{G,\pi,\zeta}(k)$, and the result follows.

When $0<\zeta<1$, we apply Proposition \ref{gantertprop} and follow
the method in \cite{Gantert}.  First, as in the
proof of Theorem \ref{marginalprop}, (\ref{dobeq}) holds, and 
$c(P^{G,\zeta}_n) =
1-C_G/n^{\zeta}$ for a constant $C_G>0$ and all $n\geq n_0(G,\zeta)$.  Then, 
$c_n = \max_{n_0(G,\zeta)\leq i\leq n} c(P_i^{G,\zeta}) =
1-C_G/n^{\zeta} <1$.  Now, $n(1-c_n) = C_G n^{1-\zeta} \uparrow
\infty$ to give the probability convergence in part (i).  For
a.s. convergence in part (ii) when $0<\zeta<1/2$, note
$$ \sum_n\frac{1}{2^n(1-c_{2^n})^2}\ =\
\sum_n\frac{1}{2^n({C_G}/{(2^n)^\zeta})^2}\ =\
\sum_n\frac{1}{C^2_G(2^{1-2\zeta})^n}\ <\ \infty . \ \ \ \  \ \ \ \ \
\ \ \ \ \ \ \ \ \ \ \ \square
$$

\sect{Proof of Theorem \ref{existencethm}.}  In this section, as $\zeta=1$ is fixed, we suppress notational dependence
on $\zeta$.  Also, as $Z_n$ takes values on the compact set
$\Delta_m$, the weak convergence in Theorem \ref{existencethm} follows by convergence
of the moments.

The next lemma establishes convergence of the first moments.

\begin{lem}
\label{1stmomlem} For $G\in\G$, $1\le k\le m$, and initial
distribution $\pi$, $$ \lim_{n\rightarrow \infty} \E_{\pi}^G\bigg(Z_{k,n}\bigg) \ =\
\nu_G(k) $$
\end{lem}

{\it Proof.} From Theorem \ref{marginalprop}, and Cesaro
convergence, $$\lim_n \E_{\pi}^G\bigg(Z_{k,n}\bigg) \ = \ \lim_n
\frac{1}{n}\sum_{i=1}^{n} \E_{\pi}^G\bigg(1_{k}(X_i)\bigg) \ = \
\lim_n \frac{1}{n}\sum_{i=1}^n \p^G_{\pi}(X_i=k) \ =\ \nu_G(k). \ \ \
\ \ \ \ \ \square
$$

We now turn to the joint moment limits in several steps,
and will assume in the following that $\gamma_1,\ldots,\gamma_m\geq
0$ with $\bar{\gamma}\ge 2$. The
first step is an ``ordering of terms.''

\begin{lem}
\label{chlemsum} For $G\in\G$, and initial distribution $\pi$, we
have
\begin{eqnarray*}
&& \lim_{n\rightarrow \infty}
\bigg|\E_{\pi}^G\bigg(Z_{1,n}^{\gamma_1}\cdots
Z_{m,n}^{\gamma_m}\bigg)\\
&&\ \ \ \ \ -\sum_{\sigma\in \S(\gamma_1,\ldots,\gamma_m)}
\frac{1}{n^{\bar{\gamma}}}\sum_{i_1=1}^{n-\bar{\gamma}+1}
\sum_{i_2>i_1}^{n-\bar{\gamma}+2}\cdots\sum_{i_{\bar{\gamma}}>
i_{\bar{\gamma}-1}}^{n} \E_{\pi}^G\bigg(\prod_{l=1}^{\bar{\gamma}}
1_{\sigma_l}(X_{i_l})\bigg)\bigg| \ = \ 0.
\end{eqnarray*}
\end{lem}

{\it Proof.} By definition of $\S(\gamma_1,\ldots,\gamma_m)$, \begin{eqnarray*}
\E_{\pi}^G\bigg( Z_{1,n}^{\gamma_1}\cdots Z_{m,n}^{\gamma_m} \bigg)
&=& \frac{1}{\bar{\gamma}!}\frac{1}{n^{\bar{\gamma}}}
\sum_{\sigma\in \S(\gamma_1,\ldots,\gamma_m) \atop 1\le
i_1,\ldots,i_{\bar{\gamma}}\le n} \E_{\pi}^G\bigg(
1_{\sigma_1}(X_{i_1})1_{\sigma_2}(X_{i_2}) \cdots
1_{\sigma_{\bar{\gamma}}}(X_{i_{\bar{\gamma}}})\bigg). \end{eqnarray*}
Note now \begin{eqnarray*} \sum_{\sigma\in \S(\gamma_1,\ldots,\gamma_m) \atop 1\le
i_1,\ldots,i_{\bar{\gamma}}\le n}1=\bar{\gamma}!n^{\bar{\gamma}} ,\
 {\rm and \ }
\sum_{\sigma\in \S(\gamma_1,\ldots,\gamma_m) \atop 1\le
i_1,\ldots,i_{\bar{\gamma}}\le n,\ {\rm
distinct}}1=\bar{\gamma}!\bar{\gamma}! {n \choose \bar{\gamma}}.
\end{eqnarray*}
Let ${\cal K}$ be those indices $\<i_1,\ldots,i_{\bar{\gamma}}\>$, $1\leq i_1,\ldots, i_{\bar{\gamma}}\leq n$
which are not distinct, that is $i_j = i_k$ for some $j\neq k$.
 Then, \begin{eqnarray*}
&&\frac{1}{\bar{\gamma}!}\frac{1}{n^{\bar{\gamma}}}\ \bigg|
\sum_{\sigma\in \S(\gamma_1,\ldots,\gamma_m) \atop 1\le
i_1,\ldots,i_{\bar{\gamma}}\le n}
\E_{\pi}^G\bigg(\prod_{l=1}^{\bar{\gamma}} 1_{\sigma_l}(X_{i_l})
\bigg) - \sum_{\sigma\in \S(\gamma_1,\ldots,\gamma_m) \atop 1\le
i_1,\ldots,i_{\bar{\gamma}}\le n,\ {\rm distinct}} \E_{\pi}^G\bigg(
\prod_{l=1}^{\bar{\gamma}}1_{\sigma_l}(X_{i_l})\bigg)\bigg |\\
&&\ \ =\ \frac{1}{\bar{\gamma}!}\frac{1}{n^{\bar{\gamma}}}
\sum_{\sigma\in \S(\gamma_1,\ldots,\gamma_m) \atop \
\<i_1,\ldots,i_{\bar{\gamma}}\>\in {\cal K}} \E_{\pi}^G\bigg(
1_{\sigma_1}(X_{i_1}) \cdots
1_{\sigma_{\bar{\gamma}}}(X_{i_{\bar{\gamma}}})\bigg) \\
&&\ \  \le\ \frac{1}{\bar{\gamma}!}
\frac{1}{n^{\bar{\gamma}}}\bigg(\bar{\gamma}!n^{\bar{\gamma}}-
\bar{\gamma}!\bar{\gamma}!{n\choose \bar{\gamma}}\bigg) \ =\ o(1).
\end{eqnarray*}
But, $$ \sum_{\sigma\in \S(\gamma_1,\ldots,\gamma_m) \atop 1\le
i_1,\ldots,i_{\bar{\gamma}}\le n,\ {\rm distinct}} \E_{\pi}^G\bigg(
\prod_{l=1}^{\bar{\gamma}}1_{\sigma_l}(X_{i_l})\bigg)\ = \
\bar{\gamma}! \sum_{\sigma\in \S(\gamma_1,\ldots,\gamma_m) \atop 1\le
i_1<\cdots<i_{\bar{\gamma}}\le n} \E_{\pi}^G\bigg(
\prod_{l=1}^{\bar{\gamma}}1_{\sigma_l}(X_{i_l})\bigg). \ \ \ \ \ \ \ \
\ \ \ \ \ \ \ \ \ 
\square $$

The next lemma replaces the initial measure with $\nu_G$.  Let
$P^G_{i,j} = \prod_{l=i}^j P^G_l$ for $1\leq i\leq j$.

\begin{lem}
\label{initdistlemcor2} For $G\in\G$ and initial distribution $\pi$,
we have
\begin{eqnarray}
\label{idlc2eq1}
&&\lim_{n\rightarrow \infty}\bigg| \sum_{\sigma\in
\S(\gamma_1,\ldots,\gamma_m)}\frac{1}{n^{\bar{\gamma}}}\sum_{i_1=1}^{n-\bar{\gamma}+1}
\sum_{i_2>i_1}^{n-\bar{\gamma}+2}\cdots\sum_{i_{\bar{\gamma}}>
i_{\bar{\gamma}-1}}^{n} \E_{\pi}^G\bigg(
\prod_{l=1}^{\bar{\gamma}}1_{\sigma_l}(X_{i_l})\bigg)\\
&&\ \ -\sum_{\sigma\in \S(\gamma_1,\ldots,\gamma_m)}
\frac{\nu_G(\sigma_1)}{n^{\bar{\gamma}}}
\sum_{i_1=1}^{n-\bar{\gamma}+1}
\sum_{i_2>i_1}^{n-\bar{\gamma}+2}\cdots\sum_{i_{\bar{\gamma}}>
i_{\bar{\gamma}-1}}^{n}
\prod_{l=1}^{\bar{\gamma}-1}P^G_{i_l+1,i_{l+1}}(\sigma_l,\sigma_{l+1})\bigg
| \ = \ 0. \nonumber 
\end{eqnarray}
\end{lem}

{\it Proof.} As $\p^{G}_\pi(X_j=t|X_i=s)=P^G_{i+1,j}(s,t)$ for
$1\leq i<j$ and $s,t\in\Sigma$, we have
\begin{eqnarray*}
&&\sum_{\sigma\in \S(\gamma_1,\ldots,\gamma_m)}
\frac{1}{n^{\bar{\gamma}}}\sum_{i_1=1}^{n-\bar{\gamma}+1}
\sum_{i_2>i_1}^{n-\bar{\gamma}+2}\cdots\sum_{i_{\bar{\gamma}}>
i_{\bar{\gamma}-1}}^{n} \E_{\pi}^G\bigg(
\prod_{l=1}^{\bar{\gamma}}1_{\sigma_l}(X_{i_l})\bigg)\\
&&\ \ \ \ = \  \sum_{\sigma\in \S(\gamma_1,\ldots,\gamma_m)}
\frac{1}{n^{\bar{\gamma}}} \sum_{i_1=1}^{n-\bar{\gamma}+1}
\sum_{i_2>i_1}^{n-\bar{\gamma}+2}\cdots\sum_{i_{\bar{\gamma}}>
i_{\bar{\gamma}-1}}^{n}\p^G_\pi(X_{i_1}=\sigma_1)\prod_{l=1}^{\bar{\gamma}-1}
P^G_{i_l+1,i_{l+1}}(\sigma_l,\sigma_{l+1})
\end{eqnarray*}
which differs from the second expression in (\ref{idlc2eq1}) by at
most \begin{eqnarray*} \sum_{\sigma\in \S(\gamma_1,\ldots,\gamma_m)}
\frac{1}{n} \sum_{i_1=1}^{n-\bar{\gamma}+1}
\bigg|\p^G_\pi(X_{i_1}=\sigma_1) -\nu_G(\sigma_1)\bigg|, \end{eqnarray*} which
vanishes by Theorem \ref{marginalprop}. \qed

We now focus on a useful class of diagonalizable matrices $$\G^* \ =
\ \bigg\{G\in \R^m\times\R^m:  {\rm \ Re} (\lambda^G_l) <1 \ {\rm for \
} 1\leq l\leq m, \ {\rm and \ } G {\rm \ is \ diagonalizable}\bigg\}$$
where $\{\lambda^G_l\}$ are the eigenvalues of $G$. As ${\rm
Re}(\lambda^G_l)\leq 0$ for $1\leq l\leq m$ when $G\in \G$, certainly all
diagonalizable $G\in \G$ belong to $\G^*$. The relevance of this
class, in the subsequent arguments, is that for $G\in \G^*$ the
resolvent $(xI-G)^{-1}$ exists for $x\geq 1$.

For $G\in\G^*$, let $V_G$ be the matrix of eigenvectors and $D_G$ be
a diagonal matrix with corresponding eigenvalue entries $D_G(i,i) =
\lambda_i^G$ so that $G = V_G D_G V_G^{-1}$. Define also for $1\leq
s,t,k\leq m$,
$$g(k;s,t)\ =\ V_G(s,k)V_G^{-1}(k,t).$$
We also denote for $a_1,\ldots,a_m\in \C$, the diagonal matrix ${\rm
Diag}(a_\cdot)$ with $i$th diagonal entry $a_i$ for $1\leq i\leq m$.
We also extend the definitions of $P_n^{G}$ and
$P_{i,j}^G$ to $G\in \G^*$ with the same formulas.  In the following, we use
the
principal value of the complex logarithm, and the usual convention
$a^{b+ic} = e^{(b+ic)\log(a)}$ for $a,b,c\in \R$ with $a>0$.


\begin{lem} For $G\in\G^*$, $s,t\in \Sigma$, and $C\leq i\leq j$ where
$C=C(G)$ is a large enough constant, 
\label{PGijlem} \begin{eqnarray*} P^G_{i,j}(s,t)&=&  \sum_{k=1}^m
\nu(k;i,j)g(k;s,t)\bigg(\frac{j}{i-1}\bigg)^{\lambda^G_k}; \end{eqnarray*}
moreover, $\nu(k;i,j)\rightarrow 1$ as $i\uparrow \infty$ uniformly over
$k$ and $j$.
\end{lem}

{\it Proof.} Straightforwardly,
$$P^G_{i,j}\ =\ V_G\prod_{k=i}^j\bigg(I+\frac{1}{k}D_G\bigg)V_G^{-1} \ = \ V_G{\rm Diag}\bigg(\prod_{k=i}^j\bigg(1+\frac{\lambda^G_\cdot}{k}\bigg)\bigg)
V_G^{-1}. $$ To expand further, we note for $z\in \C$ such that
$|z-1|<1$, we have
$$\log(z)
       \ =\ (z-1)+(z-1)^2\sum_{n=0}^\infty (-1)^{n+1}\frac{1}{n+2}(z-1)^n.
$$
and estimate $$ \bigg|\sum_{n=0}^\infty
(-1)^{n+1}\frac{1}{n+2}(z-1)^n\bigg | \ \le\ \sum_{n=0}^\infty
|z-1|^n \ =\ \bigg(1-|z-1|\bigg)^{-1}.$$ 
Let now $L$ be so large such that
$\max_{1\leq u\leq m}{|\lambda^G_u|}/{L}<1/2$.
Then, for $1\leq s\leq m$ and $k\ge L$, $$
\log\bigg(1+\frac{\lambda^G_s}{k}\bigg) \ = \
\frac{\lambda^G_s}{k}+\bigg(\frac{\lambda^G_s}{k}\bigg)^2
C_{s,k}$$ for some $C_{s,k}\in\C$ with $|C_{s,k}|\le 
(1-\max_{1\leq u\leq m}|\lambda^G_u|/L)^{-1}\leq 2$.  
Then, for $i\geq L$,
$$ \prod_{k=i}^j\bigg(1+\frac{\lambda^G_s}{k}\bigg) \ = \ \exp\bigg(
\sum_{k=i}^j\log\bigg( 1+\frac{\lambda^G_s}{k}\bigg)\bigg) \ =\
\exp\bigg(\sum_{k=i}^j\frac{\lambda^G_s}{k}+c(s;i,j)\bigg) $$ where
$c(s;i,j) = \sum_{k=i}^j({\lambda^G_s}/{k})^2 C_{s,k}$ satisfies
 $$ |c(s;i,j)|\ \leq \ 2\max_{1\leq u\leq m}|\lambda^G_u|^2
\sum_{k=i}^\infty\frac{1}{k^2}\ \to\  0  \ \ {\rm uniformly \ over \
} s \ {\rm and \ }j \ {\rm as  \ } i\uparrow \infty.$$

Let now  $$d(s;i,j)\ =\ \lambda^G_s\bigg(\sum_{k=i}^{j}
\frac{1}{k}-\int_{i-1}^{j}\frac{dx}{x}\bigg)$$ and note by the
simple estimate 
$$ \sum_{k=i}^{j}\frac{1}{k}\ < \
\int_{i-1}^{j}\frac{dx}{x}\ < \  \sum_{k=i-1}^{j-1}\frac{1}{k}$$ that
$$|d(s;i,j)| \ \leq \ \max_{1\leq u\leq m}|\lambda^G_u| \bigg(\frac{1}{j}
+ \frac{1}{i-1}\bigg) \ \leq \ \max_{1\leq u\leq
m}|\lambda^G_u|\bigg(\frac{1}{i}+\frac{1}{i-1}\bigg) \ \rightarrow \
0$$ uniformly over $j$ and $s$ as $i \uparrow\infty$. 
This allows us
to write
$$\prod_{k=i}^j\bigg(1+\frac{\lambda^G_s}{k}\bigg) \ = \
\exp\bigg(c(s;i,j)+d(s;i,j)\bigg)\bigg(\frac{j}{i-1}\bigg)^{
\lambda^G_s}.$$  Defining $\nu(s;i,j) = \exp(c(s;i,j)+d(s;i,j))$
gives after multiplying out that \begin{eqnarray*}
 P^G_{i,j} & = &
V_G{\rm Diag}\bigg(\nu(\cdot
;i,j)\bigg(\frac{j}{i-1}\bigg)^{\lambda^G_\cdot}\bigg)V_G^{-1}\\
&=& \left[ \sum_{k=1}^m \nu(k;i,j)
g(k;s,t)\bigg(\frac{j}{i-1}\bigg)^{\lambda^G_k}\right]_{s,t\in
\Sigma}\end{eqnarray*} completing the proof. \qed

To continue, define for $G\in \G^*$ the function
$T^G_{x,y}(s,t):(0,1]^2\times\Sigma^2\to \C$ by \begin{eqnarray*} T^G_{x,y}(s,t)
&=&\sum_{k=1}^m g(k;s,t)\bigg(\frac{x}{y}\bigg)^{-\lambda^G_k}. \end{eqnarray*}

\begin{lem}
\label{lem-0} For $G\in\G$,
$$
\lim_{\epsilon\downarrow 0}\lim_{n\uparrow \infty} \sum_{\sigma\in
\S(\gamma_1,\ldots,\gamma_m)}
\frac{\nu_{G}(\sigma_1)}{n^{\bar{\gamma}}}
\sum_{i_1=1}^{\gil{n\epsilon}}
\sum_{i_2>i_1}^{n-\bar{\gamma}+2}\cdots\sum_{i_{\bar{\gamma}}>
i_{\bar{\gamma}-1}}^{n}\prod_{l=1}^{\bar{\gamma}-1}
P^G_{i_l+1,i_{l+1}}(\sigma_l,\sigma_{l+1}) \ = \ 0.$$
\end{lem}

{\it Proof.} For any $\sigma\in \S(\gamma_1,\ldots,\gamma_m)$,
\begin{eqnarray*}
0&\le& \lim_\epsilon\lim_n\frac{1}{n^{\bar{\gamma}}}
\sum_{i_1=1}^{\gil{n\epsilon}}
\sum_{i_2>i_1}^{n-\bar{\gamma}+2}\cdots\sum_{i_{\bar{\gamma}}>
i_{\bar{\gamma}-1}}^{n}
\prod_{l=1}^{\bar{\gamma}-1}P^G_{i_l+1,i_{l+1}}(\sigma_l,\sigma_{l+1})\\
& \le& \lim_{\epsilon}\lim_n\frac{1}{n^{\bar{\gamma}}}(n\epsilon)
n^{\bar{\gamma}-1} \ = \ 0. \ \ \ \ \ \ \ \ \ \ \ \ \ \ \ \ \ \ \ \ \
\ \ \ \ \ \ \ \ \ \ \ \ \ \ \ \ \ \ \ \ \ \ \ \ \ \ \ \ \ \ \ \ \ \ \
\ \ \ \ \ \ \ \ \ \ \ \square\end{eqnarray*}

\begin{lem}
\label{lem2} For $G\in \G^*$, $\sigma\in
\S(\gamma_1,\ldots,\gamma_m)$, and $\epsilon>0$,
\begin{eqnarray*}&&\lim_{n\uparrow
\infty}\frac{1}{n^{\bar{\gamma}}}
\sum_{i_1=\gil{n\epsilon}+1}^{n-\bar{\gamma}+1}
\sum_{i_2>i_1}^{n-\bar{\gamma}+2}\cdots\sum_{i_{\bar{\gamma}}>
i_{\bar{\gamma}-1}}^{n}
\prod_{l=1}^{\bar{\gamma}-1}P^G_{i_l+1,i_{l+1}}(\sigma_l,\sigma_{l+1})\\
&&\ \   =\ \int_{\epsilon\leq x_1\leq x_2\leq \cdots \leq
  x_{\bar{\gamma}}\leq 1}
\prod_{l=1}^{\bar{\gamma}-1}\
 T^G_{x_l,x_{l+1}}(\sigma_l,\sigma_{l+1})\  dx_1dx_2\cdots
dx_{\bar{\gamma}}
\end{eqnarray*}
\end{lem}

{\it Proof.} From Lemma \ref{PGijlem}, as $\nu(s;i,j)\rightarrow 1$
as $i\uparrow \infty$ uniformly over $j$ and $s$, $T_{x,y}(s,t)$ is
bounded, continuous on $[\epsilon,1]^2$ for fixed $s,t$, and Riemann
convergence, we have
\begin{eqnarray*}
&&\lim_n \frac{1}{n^{\bar{\gamma}}}
\sum_{i_1=\gil{n\epsilon}+1}^{n-\bar{\gamma}+1}
\sum_{i_2>i_1}^{n-\bar{\gamma}+2}\cdots\sum_{i_{\bar{\gamma}}>
i_{\bar{\gamma}-1}}^{n} \prod_{l=1}^{\bar{\gamma}-1}
P^G_{i_l+1,i_{l+1}}(\sigma_l,\sigma_{l+1})\\
 && =  \lim_{n}
\frac{1}{n^{\bar{\gamma}}}
\sum_{i_1=\gil{n\epsilon}+1}^{n-\bar{\gamma}+1}
\cdots\sum_{i_{\bar{\gamma}}>
i_{\bar{\gamma}-1}}^{n}\prod_{l=1}^{\bar{\gamma}-1}\ \sum_{k=1}^m
\nu(k;i_l+1,i_{l+1})
g(k;\sigma_{l},\sigma_{l+1})\bigg(\frac{i_{l}/n}{i_{{l+1}}/n}\bigg)^{-\lambda^G_k}
\\
&&  = \int_{\epsilon\leq x_1\leq x_2\leq \cdots \leq
  x_{\bar{\gamma}}\leq 1}
\prod_{l=1}^{\bar{\gamma}-1}T^G_{x_l,x_{l+1}}(\sigma_l,\sigma_{l+1})
dx_1dx_2\cdots dx_{\bar{\gamma}}. \ \ \ \ \ \ \ \ \ \ \ \ \ \ \ \ \ \
\ \ \ \ \ \ \ \ \ \ \ \ \ \ \ \square
\end{eqnarray*}

\begin{lem}
\label{dct}
For $G\in \G^*$ and $\sigma \in \S(\gamma_1,\ldots,\gamma_m)$,
\begin{eqnarray*}
&&\lim_{\epsilon\downarrow 0}\int_{\epsilon\leq x_1\leq x_2\leq \cdots
  \leq x_{\bar{\gamma}}\leq 1}
\prod_{l=1}^{\bar{\gamma}-1}\
 T^G_{x_l,x_{l+1}}(\sigma_l,\sigma_{l+1})\  dx_1dx_2\cdots
dx_{\bar{\gamma}}\\
&&\ \ \ = \ \int_{0}^1\int_0^{x_{\bar{\gamma}}}\cdots\int_0^{x_2}
 T^G_{x_{\bar{\gamma}-1},x_{\bar{\gamma}}}
 (\sigma_{\bar{\gamma}-1},\sigma_{\bar{\gamma}}  ) \cdots
 T^G_{x_1,x_2}(\sigma_1,\sigma_2) dx_1dx_2\cdots dx_{\bar{\gamma}}.
\end{eqnarray*}
\end{lem}

{\it Proof.}
Let
$$f_\epsilon \ = \ 1_{\{\epsilon\leq x_1\leq x_2\leq \cdots \leq
  x_{\bar{\gamma}}\leq 1\}}\prod_{l=1}^{\bar{\gamma}-1}\
 T^G_{x_l,x_{l+1}}(\sigma_l,\sigma_{l+1}).$$
Then,
$$\lim_\epsilon f_\epsilon \ = \ 1_{\{0< x_1\leq x_2\leq \cdots \leq
  x_{\bar{\gamma}}\leq 1\}}\prod_{l=1}^{\bar{\gamma}-1}\
 T^G_{x_l,x_{l+1}}(\sigma_l,\sigma_{l+1}),$$
and $f_\epsilon$ is uniformly bounded over $\epsilon$ as
$$|f_\epsilon|
\ \leq \ \bar{f}\ = \ 1_{\{0< x_1\leq x_2\leq \cdots \leq
  x_{\bar{\gamma}}\leq 1\}}\prod_{l=1}^{\bar{\gamma}-1} \sum_{k=1}^m
|g(k;\sigma_l,\sigma_{l+1})|\bigg(\frac{x_l}{x_{l+1}}\bigg)^{-{\rm
Re}(\lambda_k^G)}.$$ The right-hand bound is integrable:  Indeed, by
Tonelli's Lemma and induction, we have
\begin{eqnarray*}
\int \bar{f} dx_1\cdots dx_{\bar{\gamma}}
&=&\int_0^1\int^{x_{\bar{\gamma}}}_0 \cdots\int^{x_2}_0
\prod_{l=1}^{\bar{\gamma}-1} \sum_{k=1}^m
|g(k;\sigma_l,\sigma_{l+1})|\bigg(\frac{x_l}{x_{l+1}}\bigg)^{-{\rm
Re}(\lambda_k^G)}
dx_1\cdots dx_{\bar{\gamma}}\\
&=& \frac{1}{\bar{\gamma}} \prod_{l=1}^{\bar{\gamma}-1} \bigg
(\sum_{k=1}^m \frac{|g(k;\sigma_l,\sigma_{l+1})|}{l-{\rm
Re}(\lambda_k^G)}\bigg).\end{eqnarray*} Hence, the lemma follows
by dominated convergence and Fubini's Theorem. \qed

\begin{lem}
\label{lem0} For $G\in\G^*$ and $\sigma\in
\S(\gamma_1,\ldots,\gamma_m)$,
$$ \int_{0}^1\int_0^{x_{\bar{\gamma}}}\cdots\int_0^{x_2}
\prod_{l=1}^{\bar{\gamma}-1}T^G_{x_l,x_{l+1}}(\sigma_l,\sigma_{l+1})
dx_1\cdots dx_{\bar{\gamma}}
 \ = \ \frac{1}{\bar{\gamma}}
\prod_{l=1}^{\bar{\gamma}-1}\bigg(lI-G\bigg)^{-1}(\sigma_{l},\sigma_{l+1}).
$$
\end{lem}

{\it Proof.} By induction, the integral equals
\begin{eqnarray*}
&&\int_{0}^1\int_0^{x_{\bar{\gamma}}}\cdots\int_0^{x_2}
T^G_{x_{\bar{\gamma}-1},x_{\bar{\gamma}}}
(\sigma_{\bar{\gamma}-1}, \sigma_{\bar{\gamma}}  ) \cdots
T^G_{x_1,x_2}(\sigma_1,\sigma_2) dx_1\cdots dx_{\bar{\gamma}}\\
&&\ \ \ \ \ \ \ \ \ \ \ \ \ \ \ \ \ \ \ \ \ \ \ \ = \
\frac{1}{\bar{\gamma}} \prod_{l=1}^{\bar{\gamma}-1}
\bigg(\sum_{k=1}^m
\frac{g(k;\sigma_{l},\sigma_{l+1})}{l-\lambda^G_k}\bigg).
\label{lem0eq1}
\end{eqnarray*}
However, for $x\ge 1$, we have
$$
\bigg(xI-G\bigg)^{-1}(s,t) \ =\ V_G\bigg(xI-D_G\bigg)^{-1}V_G^{-1}(s,t)\ =\ \sum_{k=1}^m
\frac{g(k;s,t)}{x-\lambda^G_k}$$
to finish the identification. \qed

At this point, by straightforwardly combining the previous lemmas,
we have proved Theorem \ref{theorem1} for $G\in \G$ diagonalizable.
The method in extending to non-diagonalizable generators is
accomplished by
approximating with suitable ``lower'' and ``upper'' diagonal
matrices.

\begin{lem}
\label{lem1} For $G\in\G$, \begin{eqnarray} &&\lim_{n\rightarrow
\infty} \sum_{\sigma\in \S(\gamma_1,\ldots,\gamma_m)}
\frac{\nu_G(\sigma_1)}{n^{\bar{\gamma}}}
\sum_{i_1=1}^{n-\bar{\gamma}+1}
\sum_{i_2>i_1}^{n-\bar{\gamma}+2}\cdots\sum_{i_{\bar{\gamma}}>
i_{\bar{\gamma}-1}}^{n}
\prod_{l=1}^{\bar{\gamma}-1}P^G_{i_l+1,i_{l+1}}(\sigma_l,\sigma_{l+1})\nonumber\\
&&\ \ \ \ \ \ \ \ \ \ \ \ \ \ \ \ \ \ \ \ \ \  =\ \frac{1}{\bar{\gamma}}\sum_{\sigma\in
\S(\gamma_1,\ldots,\gamma_m)}  \nu_G(\sigma_1)
\prod_{l=1}^{\bar{\gamma}-1} \bigg(lI-G\bigg)^{-1}(\sigma_{l},\sigma_{l+1}).
\label{lem1eq1}\end{eqnarray}
\end{lem}

{\it Proof.} For an $m\times m$ matrix $A$, let $G[A]=G+A$. Let
$\|\cdot\|_{\rm M}$ be the matrix norm $\|A\|_{\rm M} =
\max\{|A(s,t)|: 1\leq s,t\leq m\}$.  Now, for small $\epsilon>0$, choose
matrices $A_1$ and $A_2$ with non-negative entries so that
$\|A_1\|_{\rm M},\|A_2\|_{\rm M}<\epsilon$, $I+G[-A_1]/l,I+G[A_2]/l$ have
positive entries for all $l$ large enough, and
$G[-A_1],G[A_2]\in\G^*$: This last condition can be met as (1) the spectrum
varies continuously with respect to the matrix norm $\|\cdot\|_{\rm
M}$ (cf. Appendix D \cite{HJ}), and (2) diagonalizable real matrices
are dense (cf. Theorem 1 \cite{Hartfiel}).   

Then, for $s,t\in \Sigma$, and $l$ large enough, we
have
$0<(I+G[-A_1]/l)(s,t) \leq (I+G/l)(s,t)\leq (I+G[A_2]/l)(s,t)$.  Hence,
for $i\leq j$ with $i$ large enough,
$$ P^{G[-A_1]}_{i,j}(s,t)
\ \le\ P^{G}_{i,j}(s,t) \ \le\
P^{G[A_2]}_{i,j}(s,t). $$ By Lemmas
\ref{lem-0}, \ref{lem2}, \ref{dct} and \ref{lem0}, the left-side of
(\ref{lem1eq1}), that is in terms of liminf and limsup, is bounded below and
above by
$$ \sum_{ \sigma\in \S(\gamma_1,\ldots,\gamma_m)}
\frac{1}{\bar{\gamma}}
\nu_{G}(\sigma_1)\prod_{l=1}^{\bar{\gamma}-1}\bigg(
lI-G[-A_1]\bigg)^{-1}(\sigma_{l},\sigma_{l+1}),$$ and $$ \sum_{ \sigma\in
\S(\gamma_1,\ldots,\gamma_m)} \frac{1}{\bar{\gamma}}
\nu_{G}(\sigma_1)\prod_{l=1}^{\bar{\gamma}-1}
\bigg(lI-G[A_2]\bigg)^{-1}(\sigma_{l},\sigma_{l+1})$$ respectively.  On the
other hand, for $\sigma\in \S(\gamma_1,\ldots,\gamma_m)$, both
$$
\prod_{l=1}^{\bar{\gamma}-1}
\big(lI-G[-A_1]\big)^{-1}(\sigma_{l},\sigma_{l+1}),
\prod_{l=1}^{\bar{\gamma}-1} \big(lI-G[A_2]\big)^{-1}(\sigma_l,\sigma_{l+1})
\ \rightarrow \
\prod_{l=1}^{\bar{\gamma}-1}\big(lI-G\big)^{-1}(\sigma_{l},\sigma_{l+1})$$
 as $\epsilon\to 0$, completing the proof. \qed

\sect{Proof of Theorem \ref{Dirichletthm}}

The proof follows by evaluating the moment expressions in Theorem
\ref{theorem1} when $G = \Theta$ as those corresponding to the
Dirichlet distribution with parameters $\theta_1,\ldots,\theta_m$ (\ref{dirichletmoments}).

%
%
%
%
%

\begin{lem}
\label{4.1}
The stationary distribution $\nu_\Theta$ is given by $\nu_\Theta(l)
= \theta_l/\bar{\theta}$ for $l\in \Sigma$.

Also, for $2\leq l\leq \bar{\gamma}$, let $F_l$ be the $m\times m$
matrix with entries \begin{eqnarray*} {F}_l(j,k)&=&\left\{
\begin{array}{rl}
\theta_k& \ {\rm for \ }k\not=j\\
\theta_{j}+l-1&\ {\rm for \ }k=j.  \end{array}\right. \end{eqnarray*} Then,
$$\bigg(lI-\Theta\bigg)^{-1} \ =\  \frac{1}{l(l+\bar{\theta})}F_{l+1}.$$
\end{lem}

{\it Proof.} The form of $\nu_\Theta$ follows by inspection. For
the second statement, write $F_{l+1} = lI+\hat{\Theta}$ where the
matrix $\hat{\Theta}$ has $i$th column equal to
$\theta_i(1,\ldots,1)^t$.  Then, also $\Theta = \hat{\Theta}-\bar{\theta}I$.  As $(1,\ldots, 1)^t$ is an eigenvector of $\Theta$
with eigenvalue $0$, we see $(lI-\Theta)(lI+\hat{\Theta}) =
(l^2+l\bar{\theta})I$ finishing the proof. \qed

The next statement is an immediate corollary of Theorem
\ref{existencethm} and Lemma \ref{4.1}.
\begin{lem}
\label{Dirlem0} The $\mu_\Theta$-moments satisfy
$E_{\mu_\Theta}[x_i] = \theta_i/\bar{\theta}$ for $1\leq i\leq m$ and,
when $\bar{\gamma}\geq 2$,
 \begin{eqnarray*}
E_{\mu_\Theta}\bigg[\prod_{i=1}^m x_i^{\gamma_i}\bigg] & = &
\sum_{\sigma\in \S(\gamma_1,\ldots,\gamma_m)} \nu_\Theta(\sigma_1)
\frac{1}{\bar{\gamma}}
\prod_{l=1}^{\bar{\gamma}-1}\bigg(lI-\Theta\bigg)^{-1}(\sigma_{l},\sigma_{l+1})\\
&=& \sum_{\sigma\in \S(\gamma_1,\ldots,\gamma_m)}
\frac{\theta_{\sigma_1}\prod_{l=2}^{\bar{\gamma}}
{F}_l(\sigma_{l-1},\sigma_{l})}
{\bar{\gamma}!\prod_{l=0}^{\bar{\gamma}-1}(\bar{\theta}+l)}. \end{eqnarray*}
\end{lem}


We now evaluate the last expression of Lemma \ref{Dirlem0} by first
specifying of the value of $\sigma_{\bar{\gamma}}$.  Recall, by
convention $\theta_l\cdots(\theta_l+\gamma_l-1)=1$ when $\gamma_l=0$
for $1\le l\le m$.

\begin{lem}
\label{lem2base} For $\bar{\gamma}\geq 2$ and $1\le k\le m$,
\begin{equation}
\label{lem2baseeq1}\sum_{\sigma\in \S(\gamma_1,\ldots,\gamma_m)\atop
\sigma_{\bar{\gamma}}=k} \theta_{\sigma_{1}}
\prod_{l=2}^{\bar{\gamma}} {F}_l(\sigma_{l-1},\sigma_{l})  \ =\
\gamma_k(\bar{\gamma}-1)!\prod_{l=1}^m\theta_l\cdots(\theta_l+\gamma_l-1).
\end{equation}
\end{lem}

{\it Proof.} The proof will be by induction on $\bar{\gamma}$.
\vskip .1cm

\noindent {\it Base Step: $\bar{\gamma}= 2$.}  If $\gamma_k=1$ and
$\gamma_i=1$ for $i\neq k$, the left and right-sides of
(\ref{lem2baseeq1}) both equal
$\theta_i{F}_2(i,k)=\theta_i\theta_k$. If $\gamma_k=2$, then the
left and right-sides of (\ref{lem2baseeq1}) equal
$2\theta_k{F}_2(k,k)=2\theta_k(\theta_k+1)$. \vskip .1cm

\noindent {\it Induction Step.} Without loss of generality and to
ease notation, let $k=1$. Then, by specifying the next-to-last
element $\sigma_{\bar{\gamma}-1}$, and simple counting, we have
\begin{eqnarray*}
\sum_{\sigma\in \S(\gamma_1,\ldots,\gamma_m)\atop
\sigma_{\bar{\gamma}}=1}
\theta_{\sigma_{1}}\prod_{l=2}^{\bar{\gamma}}
{F}_l(\sigma_{l},\sigma_{l-1})
&=&\gamma_1(\theta_1+\bar{\gamma}-1)\sum_{\sigma\in
\S(\gamma_1-1,\ldots,\gamma_m) \atop \sigma_{\bar{\gamma}-1}=1}
\theta_{\sigma_1}\prod_{l=2}^{\bar{\gamma}-1}
{F}_l(\sigma_{l},\sigma_{l-1})\\
&&\ \ \ +\sum_{j=2}^m \gamma_1\theta_1 \sum_{\sigma\in
\S(\gamma_1-1,\ldots,\gamma_m)\atop \sigma_{\bar{\gamma}-1}=j}
\theta_{\sigma_{1}}\prod_{l=2}^{\bar{\gamma}-1}
{F}_l(\sigma_{l},\sigma_{l-1}).
\end{eqnarray*}
We now use induction to evaluate the right-side above as
\begin{eqnarray*}
&&\theta_1\cdots(\theta_1+\gamma_1-2)\prod_{i=2}^m\theta_i\cdots(\theta_i+\gamma_i-1)\\
&&\ \ \ \ \ \ \ \ \ \ \ \ \ \ \ \ \times
\bigg\{\gamma_1(\theta_1+\bar{\gamma}-1)(\gamma_1-1)(\bar{\gamma}-2)!
+\sum_{j=2}^m \gamma_1\theta_1\gamma_j(\bar{\gamma}-2)!\bigg\}\\
&&\ \ \ \ \ \ \ =\
\theta_1\cdots(\theta_1+\gamma_1-2)\prod_{i=2}^m\theta_i\cdots(\theta_i+\gamma_i-1)\\
&&\ \ \ \ \ \ \ \ \ \ \ \ \ \ \ \ \times
\bigg\{\gamma_1(\theta_1+\bar{\gamma}-1)(\gamma_1-1)(\bar{\gamma}-2)!
+\gamma_1\theta_1(\bar{\gamma}-\gamma_1)(\bar{\gamma}-2)!\bigg\}\\
&&\ \ \ \ \ \ \ =\
\theta_1\cdots(\theta_1+\gamma_1-2)\prod_{i=2}^m\theta_i\cdots(\theta_i+\gamma_i-1)\\
&&\ \ \ \ \ \ \ \ \ \ \ \ \ \ \ \ \times
\gamma_1(\bar{\gamma}-2)!\bigg\{(\theta_1+\gamma_1-1)(\bar{\gamma}-1)\bigg\}\\
&&\ \ \ \ \ \ \ =\
\gamma_1(\bar{\gamma}-1)!\prod_{l=1}^m\theta_l\cdots(\theta_l+\gamma_l-1).
\ \ \ \ \ \ \ \ \ \ \ \ \ \ \ \ \ \ \ \ \ \ \ \ \ \ \ \ \ \ \ \ \ \ \
\ \ \ \ \ \ \ \ \ \ \ \ \ \ \ \ \ \square
\end{eqnarray*}

By now adding over $1\le k\le m$ in the previous lemma, we finish
the proof of Theorem \ref{Dirichletthm}.
\begin{lem}
\label{lem2count} When $\bar{\gamma}\geq 2$, \begin{eqnarray*} \sum_{\sigma\in
\S(\gamma_1,\ldots,\gamma_m)} 
\frac{
\theta_{\sigma_{1}} \prod_{l=2}^{\bar{\gamma}}
{F}_l(\sigma_{l-1},\sigma_{l})}
{\bar{\gamma}!\prod_{l=0}^{\bar{\gamma}-1}(\bar{\theta}+l)} 
&=& \frac{
\prod_{l=1}^m\theta_l\cdots(\theta_l+\gamma_l-1)}
{\prod_{l=0}^{\bar{\gamma}-1}(\bar{\theta}+l)}. \end{eqnarray*}
\end{lem}

{\it Proof.}  \begin{eqnarray*} \sum_{\sigma\in \S(\gamma_1,\ldots,\gamma_m)}
\frac{ \theta_{\sigma_{1}}
\prod_{l=2}^{\bar{\gamma}} {F}_l(\sigma_{l-1},\sigma_{l})
}{\bar{\gamma}!\prod_{l=0}^{\bar{\gamma}-1}(\bar{\theta}+l)} &=&
\sum_{k=1}^m\sum_{\sigma\in \S(\gamma_1,\ldots,\gamma_m) \atop
\sigma_{\bar{\gamma}}=k} 
\frac{
\theta_{\sigma_{1}} \prod_{l=2}^{\bar{\gamma}}
{F}_l(\sigma_{l-1},\sigma_{l})
}{\bar{\gamma}!\prod_{l=0}^{\bar{\gamma}-1}(\bar{\theta}+l)}
\\
&=&\frac{\sum_{k=1}^m\gamma_k(\bar{\gamma}-1)!} {\bar{\gamma}!}
\frac{\prod_{l=1}^m\theta_l\cdots(\theta_l+\gamma_l-1) }{
\prod_{l=0}^{\bar{\gamma}-1}(\bar{\theta}+l)
}\\
&=& \frac{\prod_{l=1}^m\theta_l\cdots(\theta_l+\gamma_l-1)}
{\prod_{l=0}^{\bar{\gamma}-1}(\bar{\theta}+l)}. \ \ \ \ \ \ \ \ \ \ \
\ \ \ \ \ \ \ \ \ \ \ \ \ \ \square \end{eqnarray*}  

\sect{Proof of Theorem \ref{structurethm} (1)}
Let ${\bf p}=\<p_1,\ldots,p_m\>\in{\rm Int}\Delta_m$ be a point in the
simplex with $p_i>0$ for $1\leq i\leq m$.  For $\epsilon>0$ small,
let $B({\bf p},\epsilon)\subset {\rm Int}\Delta_m$ be a ball with radius
$\epsilon$ and center ${\bf p}$.  To prove Theorem \ref{structurethm} (1),
it is enough to show for all large $n$ the lower bound
$$\p_\pi^G\bigg({\bf Z}_n \in B({\bf p},\epsilon)\bigg)\ > \ C({\bf p},\epsilon)\ >\
0.$$

To this end, let $\bar{p}_0=0$ and $\bar{p}_i = \sum_{l=1}^i p_l$ for $1\leq
i\leq m$.  Also, define, for $1\leq k\leq l$, ${\bf X}_k^l=\<X_k,\ldots,
X_l\>$.
Then, there
exist small $\delta,\beta>0$
 such that
\begin{eqnarray}
\label{keystep}
&&\big\{{\bf Z}_n\in B({\bf p},\epsilon)\big\}\\
&&\ \ \supset \
\cup_{0\leq k_1,\ldots,k_m\leq \gil{n\beta}} \bigg\{ \bigg\{
{\bf X}_{\gil{n\delta}}^{\gil{n\bar{p}_1}-k_1} =\vec{1}\bigg\}\cap\bigg(
\cap_{j=2}^{m} \bigg\{
{\bf X}_{\gil{n\bar{p}_{j-1}}-\bar{k}_{j-1}+1}^{\gil{n\bar{p}_j}-\bar{k}_j}
=\vec{j}\bigg\} \bigg)\bigg\}\nonumber
\end{eqnarray}
where $\bar{k}_a = \sum_{l=1}^a k_l$, and $\vec{i}$ is a vector with
all coordinates equal to $i$ of the appropriate length.  The last
event represents the process being in the fixed location $j$ for
times $\gil{n\bar{p}_{j-1}} - \bar{k}_{j-1} +1$ to $\gil{n\bar{p}_j}
- \bar{k}_j$ for $1\leq j\leq m$ where we take $1-\bar{k}_0 =
\gil{n\delta}$.

Now, as $G$ has strictly negative diagonal entries, $C_1=\max_{s}|G(s,s)|>0$, and so for
all large $n$,
$$
\p^G_\pi\bigg({\bf X}_{\gil{n\bar{p}_{i-1}}-\bar{k}_{i-1}+1}^{\gil{n\bar{p}_{i}}-\bar{k}_i}
=\vec{i}\big| {\bf X}_{\gil{n\bar{p}_{i-1}}-\bar{k}_{i-1}+1}=i\bigg) \
\ge\ \prod_{j=\gil{n\delta}}^n 1-\frac{C_1}{j}\
 \geq \ \frac{\delta^{C_1}}{2}.$$
Also, as $G$ has positive nondiagonal entries, $C_2=\min_{s} G(s,s+1)>0$.  Then, \begin{eqnarray*}
\p^G_\pi\bigg( X_{\gil{n\bar{p}_{i-1}}-\bar{k}_{i-1}+1}=i\big| X_{\gil{n\bar{p}_{i-1}}-\bar{k}_{i-1}}=i-1\bigg)
&\ge& \frac{C_2}{\gil{n\bar{p}_{i-1}}-\bar{k}_{i-i}+1}. \end{eqnarray*}

Hence, for all large $n$, as $\p_\pi^G(X_{\gil{n\delta}} = 1) \geq
\nu_G(1)/2$ (Theorem \ref{marginalprop}), \begin{eqnarray*} &&\p^G_\pi\bigg(
{\bf Z}_n\in B({\bf p},\epsilon)\bigg)\\
&&\ \ \ \ \ge\ \sum_{0\leq k_1,\ldots,k_m\leq \gil{n\beta}}
\p^G_\pi\bigg( \bigg\{ {\bf X}_{\gil{n\delta}}^{\gil{np_{1}}-k_1}
=\vec{1}\bigg\}\cap\bigg( \cap_{j=2}^{m} \bigg\{
{\bf X}_{\gil{n\bar{p}_{j-1}}-\bar{k}_{j-1}+1}^{\gil{n\bar{p}_{j}}-\bar{k}_{j}}
=\vec{j}
\bigg\}\bigg)\bigg)\\
&&\ \ \ \ \ge\ \bigg[\frac{\delta^{C_1}}{2}\bigg]^m \sum_{0\leq
k_1,\ldots,k_m\leq \gil{n\beta}}\frac{\nu_G(1)}{2}\prod_{j=2}^{m}
\frac{C_2}{\gil{n\bar{p}_{j-1}}-\bar{k}_{j-1} +1}\\
&&\ \ \ \ \ge\ \bigg[\frac{\delta^{C_1}}{2}\bigg]^m \sum_{0\leq
k_1,\ldots,k_m\leq \gil{n\beta}}\frac{\nu_G(1)}{2}\prod_{j=2}^{m}
\frac{C_2}{\gil{n\bar{p}_{j-1}}-k_{j-1} +1}\\
&&\ \ \ \ \ge\
\frac{\nu_G(1)}{4}\bigg[\frac{C_2\delta^{C_1}}{2}\bigg]^m
\prod_{j=2}^{m}
\log\bigg(\frac{\gil{n\bar{p}_{j-1}}}{\gil{n\bar{p}_{j-1}}-\gil{n\beta}}\bigg)\\
&&\ \ \ \ \geq\
\frac{\nu_G(1)}{8}\bigg[\frac{C_2\delta^{C_1}}{2}\bigg]^m
\prod_{j=2}^{m}\log\bigg(\frac{\bar{p}_{j-1}}{\bar{p}_{j-1}-\beta}\bigg).
\ \ \ \ \ \ \ \ \ \ \ \ \ \ \ \ \ \ \ \ \ \ \ \ \ \ \ \ \ \ \ \ \ \ \
\ \ \ \  \ \ \ \ \ \ \ \ \ \square
\end{eqnarray*}


\sect{Proof of Theorem \ref{structurethm} (2)}
The proof of Theorem \ref{structurethm} (2) follows from the next
two propositions.
\begin{prop} For $G\in \G$, the $m$ vertices of
$\Delta_m$, ${\bf 1},\ldots,{\bf m}$, are not atoms.
\end{prop}

{\it Proof.} From Theorem \ref{existencethm}, moments $\alpha_{l,k}
= E_{\mu_G}[(x_l)^k]$ satisfy $\alpha_{l,k+1} = (I-G/k
)^{-1}(l,l)\alpha_{l,k}$ for $1\leq l\leq m$ and $k\geq 1$. By the
inverse adjoint formula, for large $k$,
$$\bigg(I-G/k\bigg)^{-1}(l,l) \ = \ \frac{1-
\frac{1}{k}({\rm Tr}(G) - G(l,l))}{1-{\rm Tr}(G)/k} +O(k^{-2}) \ = \
1+\frac{G(l,l)}{k} +O(k^{-2}).$$  As $G\in \G$, $G(l,l)<0$.  Hence, $\alpha_{l,k}$ vanishes at polynomial
rate $\alpha_{l,k} \sim k^{G(l,l)}$.  In particular, as
$\mu_G(\{{\bf l}\})\leq E_{\mu_G}[(x_l)^k]\rightarrow 0$ as
$k\rightarrow \infty$, the point ${\bf l}$ cannot be an atom of the
limit distribution. \qed

Fix for the remainder ${\bf p}\in\Delta_m\setminus\{{\bf
1},\ldots,{\bf m}\}$, and define $\check{p} = \min\{p_i: p_i> 0,
1\leq i\leq m\}>0$. Let also $0<\delta< \check{p}/2$, and consider
$B({\bf p},\delta)= \{{\bf x}\in\Delta_m: |{\bf p} -{\bf x}|<\delta\}$.
\begin{prop}
\label{absprop} For $G\in \G$, there is a constant $C=C(G,{\bf p},m)$ such
that
$$\mu_G\bigg(B({\bf p},\delta)\bigg) \ \leq \ C \ \log \bigg(\frac{\check{p} + 2\delta }
{\check{p} - \delta }\bigg). $$
\end{prop}

Before proving Proposition \ref{absprop}, we will need some notation
and lemmas. We will say a ``switch'' occurs at time $1< k\leq n$ in
the sequence $\omega^n =\<\omega_{1},\ldots,\omega_n\>\in \Sigma^n$
if $\omega_{k-1} \neq \omega_{k}$. For $0\leq j\leq n-1$, let
$$T(j) \ =\
\bigg\{\omega^n: \omega^n \ {\rm has \ exactly \ } j {\rm \
switches}\bigg\}.$$ 
Note as ${\bf p}\in \Delta_m\setminus\{{\bf
1},\ldots,{\bf m}\}$ at least two coordinates of ${\bf p}$ are positive.
Then, as $\delta <\check{p}/2$, when $(1/n)\sum_{i=1}^n
\<1_{1}(\omega_i),\ldots,1_{m}(\omega_i)\> \in B({\bf p},\delta)$, at least
one switch is in $\omega^n$.

For $j\geq 1$ and a path in $T(j)$, let $\alpha_1,\ldots,\alpha_{j}$
denote the $j$ switch times in the sequence; let also
$\theta_1,\ldots,\theta_{j+1}$ be the $j+1$ locations visited by the
sequence. We now partition $\{\omega^n: (1/n)\sum_{i=1}^n
\<1_1(\omega_i),\ldots,1_m(\omega_i)\>\in B({\bf p},\delta)\}\cap T(j)$ into non-empty sets $A_j({\bf
U},{\bf V})$ where ${\bf U} = \<U_1,\ldots,U_{j-1}\>$ and ${\bf V}=
\<V_1,\ldots,V_{j+1}\>$ denote possible switch times (up to the
$j-1$st switch time) and visit locations respectively:
\begin{eqnarray*}
A_j({\bf U},{\bf V})  &=&  \bigg\{\omega^n: \omega^n\in T(j),
\frac{1}{n}\sum_{i=1}^n \<1_1(\omega_i),\ldots,1_m(\omega_i)\>\in
B({\bf p},\delta),\\
&&\ \ \ \ \ \ \ \ \ \ \ \ \ \ \ \ \ \ \ \ \ \ \ \   \alpha_i = U_i,
\theta_k = V_k \ {\rm for \ }1\leq i \leq j-1, 1\leq k\leq
j+1\bigg\}.\end{eqnarray*} In this decomposition, paths in $A_j({\bf
U}, {\bf V})$ are in $1:1$ correspondence with $j$th switch times
$\alpha_j$--the only feature allowed to vary.

Now, for each set $A_j({\bf U},{\bf V})$, we define a path
$\eta(j,{\bf U}, {\bf V}) = \<\eta_{1},\ldots,\eta_n\>$ where the
last $j$th switch is ``removed,''
$$ \eta_l = \ \left\{ \begin{array}{ll}
V_1 & {\rm for \ }1 \leq l< U_1\\
V_k & {\rm for \ }U_{k-1} \leq l< U_k, 2\leq k \leq j-1\\
V_j & {\rm for \ }U_{j-1}\leq l\le n. \end{array} \right. $$  Note
that the sequence $\eta(j,{\bf U},{\bf V})$ belongs to $T(j-1)$, can
be obtained no matter the location $V_{j+1}$ (which could range on
the $m$ values in the state space), and is in $1:1$ correspondence
with pair $\<U_1,\ldots,U_{j-1}\>$ and $\<V_1,\ldots,V_j\>$. In
particular, recalling ${\bf
X}^n_1 = \<X_1,\ldots,X_n\>$ denotes the coordinate sequence up to time $n$, we have
\begin{equation}
\label{doubling} \sum_{{\bf U,V}} \p^G_\pi\bigg({\bf X}^n_1=\eta(j,{\bf
U,V})\bigg) \ \leq \ m \ \p^G_\pi\bigg({\bf X}^n_1\in T(j-1)\bigg)
\end{equation}
where the sum is over all ${\bf U},{\bf V}$ corresponding to the
decomposition into sets $A_j({\bf U,V})$ of
 $\{\omega^n: (1/n)\sum_{i=1}^n
\<1_1(\omega_i),\ldots,1_m(\omega_i)\>\in B({\bf p},\delta)\}\cap T(j)$.

The next lemma estimates the location of the last switch time
$\alpha_j$, and the size of the set $A_j({\bf U,V})$.  The proof is
deferred to the end.
\begin{lem}
\label{lastswitchtimelemma} On $A_j({\bf U},{\bf V})$, we have
$\lceil n(\check{p}-\delta)+1\rceil  \leq  \alpha_j$.  Also,
$|A_j({\bf U},{\bf V})|  \leq  \lfloor 2n\delta +1\rfloor$.
\end{lem}
A consequence of these bounds on the position and cardinality of
$\alpha_j$'s associated to a fixed set $A_j({\bf U,V})$, is that
\begin{equation}\label{lastswitchtimebound} \sum' \frac{1}{U_j} \
 \leq\ \sum_{k=\lceil n(\check{p}- \delta)+1\rceil}^{\lceil n(\check{p}+\delta)+2\rceil }\frac{1}{k} \ \leq\
\log \bigg(\frac{\check{p}+ \delta +3/n}{\check{p}-\delta}\bigg)
\end{equation}
where $\sum'$ refers to adding over all last switch times $U_j$
associated to paths in $A_j({\bf U},{\bf V})$.

Let now $\hat{G} = \max\{|G(i,j)|: 1\leq i,j\leq m \}$.

\begin{lem}
\label{etalem} For $\omega^n\in A_j({\bf U},{\bf V})$ such that
$\alpha_j = U_j$, and all large $n$, we have
\begin{equation}
\label{equation1} \p^G_\pi\bigg({\bf X}^n_1=\omega^n\bigg) \ \le \
\frac{\hat{G}(\check{p}/2)^{-2\hat{G}}}{U_j}\ \p^G_\pi\bigg({\bf
X}^n=\eta(j,{\bf U},{\bf V})\bigg).
\end{equation}
\end{lem}

{\it Proof.}  The path $\eta(j,{\bf U,V})$ differs from $\omega^n$
only in that there is no switch at time $U_j$. Hence,
$$\frac{\p^G_\pi\big({\bf X}^n=\omega^n\big)}{\p^G_\pi\big({\bf X}^n=\eta(j,{\bf U},{\bf
V})\big)} \ = \ \frac{G(V_j,V_{j+1})}{U_j(1+G(V_j,V_j)/U_j)}
\prod_{l=U_j +1}^{n}
\bigg(\frac{1+G(V_{j+1},V_{j+1})/l}{1+G(V_j,V_j)/l}\bigg).$$ 
Now
bounding $G(V_j,V_{j+1})\leq \hat{G}$, $1+G(V_{j+1},V_{j+1})/l
\leq 1$, $1+G(V_j,V_j)/l\geq 1-\hat{G}/l$, and noting $U_j \geq
n(\check{p}-\delta) +1$ (by Lemma \ref{lastswitchtimelemma}),
$-\ln(1-x)\leq 2x$ for $x>0$ small, and $\delta<\check{p}/2$, give for
large $n$, 
$$
\frac{G(V_j,V_{j+1})}{1+G(V_j,V_j)/U_j} \prod_{l=U_j +1}^{n}
\bigg(\frac{1+G(V_{j+1},V_{j+1})/l}{1+G(V_j,V_j)/l}\bigg) \leq
{\hat{G}} \bigg(\frac{n}{n(\check{p}-\delta)}\bigg)^{2\hat{G}} \leq  
\hat{G}(\check{p}/2)^{-2\hat{G}}.\ \ \ \   \square $$

 {\it Proof
of Proposition \ref{absprop}.} By decomposing over number of
switches $j$ and on
 the structure of the paths with $j$ switches, estimates
 (\ref{equation1}), (\ref{lastswitchtimebound}), comment (\ref{doubling}),
 and $\sum_j \p^G_\pi \big({\bf X}^n\in T(j-1)\big) \leq 1$, we have
 for all large $n$,
\begin{eqnarray*}
\p^G_\pi\bigg({\bf Z}_n\in B({\bf p}, \delta)\bigg)&=& \sum_{j=1}^{n-1}
\p^G_\pi\bigg({\bf Z}_n\in B({\bf p}, \delta), {\bf X}^n\in T(j)\bigg)\\
&=& \sum_{j=1}^{n-1}
\sum_{{\bf U,V}}\p^G_\pi\bigg(A_j({\bf U},{\bf V})\bigg)\\
&\leq& \sum_{j=1}^{n-1} \sum_{{\bf U,V}}\sum'\
\frac{C(G,{\bf p})}{U_j}\ \p^G_\pi\bigg({\bf X}^n=\eta(j,{\bf U},{\bf V})\bigg)\\
&\leq& C(G,{\bf p}) \log \bigg(\frac{\check{p}+2\delta} {\check{p}-\delta}
\bigg)\sum_{j=1}^{n-1} \sum_{{\bf U,V}}\ \p^G_\pi\bigg({\bf
X}^n=\eta(j,{\bf U},{\bf V})\bigg) \\
&\leq& m C(G,{\bf p}) \log \bigg(\frac{\check{p}+2\delta}
{\check{p}-\delta} \bigg)\sum_{j=1}^{n-1}\ \p^G_\pi\bigg({\bf X}^n\in T(j-1)\bigg)
\\
&\leq& C(G,{\bf p},m) \log \bigg(\frac{\check{p}+2\delta}
{\check{p}-\delta} \bigg).
\end{eqnarray*}
The proposition follows by taking limit on $n$, and weak convergence.\qed

{\it Proof of Lemma \ref{lastswitchtimelemma}.} For a path $\omega^n
\in A_j({\bf U,V})$ and $1\le k\le j+1$, let $\tau_k$ be the number
of visits to state $V_k$ (some $\tau_k$'s may be the same if $V_k$ is repeated).  For $1\leq i\leq \tau_k$, let
$\underline{n}^k_i$ and $\overline{n}^k_i$ be the start and end of
the $i$th visit to $V_k$. Certainly, $ \sum_{i=1}^n
1_{V_k}(\omega_i)  = \sum_{i=1}^{\tau_k}
(\overline{n}^k_{i}-\underline{n}^k_i +1)$. Moreover, as
$({1}/{n})\sum_{i=1}^n \<1_1(\omega_i),\ldots,1_m(\omega_i)\>\in
B({\bf p},\delta)$, we have $ |(1/n)\sum_{i=1}^{n}
1_{V_k}(\omega_i)-p_{V_k}| \le \delta $, and so
\begin{equation}
\label{start-end-bounds} n(p_{V_k}-\delta)\ \le \
\sum_{i=1}^{\tau_k} \big(\overline{n}^k_{i}-\underline{n}^k_i
+1\big) \ \le \ n(p_{V_k}+\delta). \end{equation} Hence, as 
the disjoint sojourns
$\{[\underline{n}_i^k,\overline{n}_i^k]: 1\leq i\leq \tau_k\}$ occur
between times $1$ and $\overline{n}^k_{\tau_k}$, their total sum
length is less than $\overline{n}^k_{\tau_k}$, and we deduce
$n(p_{V_k}-\delta)\leq \overline{n}^k_{\tau_k}$.

Now, for ${\bf p}\in \Delta_m\setminus\{{\bf 1},\ldots,{\bf m}\}$, at
least one of the $\{p_{V_i}:V_i\neq V_{j+1}, 1\le i\le j\}$ is positive: Indeed,
there are two coordinates of ${\bf p}$, say $p_s$ and $p_t$, which are
positive.  Say $V_{j+1}\neq s$; then, as $(1/n)\sum_{i=1}^n
1_s(\omega_i) = (1/n)\sum_{i=1}^{\alpha_j-1} 1_s(\omega_i)$,
$|(1/n)\sum_{i=1}^{\alpha_j-1} 1_s(\omega_i) - p_s|\leq \delta$, and
$p_s-\delta>0$, the path must visit state $s$ before time
$\alpha_j$, e.g. $V_i=s$ for some $1\leq i\leq j$. 

Then,
from the deduction just after (\ref{start-end-bounds}), we have
$$n(\check{p}-\delta) \ \leq \ n\max_{V_i\neq V_{j+1}\atop 1\leq i\leq
  j}(p_{V_i}-\delta)\ \leq \ \max_{V_i\neq V_{j+1}\atop 1\leq i\leq j} \overline{n}^i_{\tau_i}\
\leq \ \overline{n}^j_{\tau_j} \ =\  \alpha_j -1$$ giving the first
statement.

For the second statement, note that $-\underline{n}_{\tau_j}^j +
\sum_{i=1}^{\tau_j -1}(\overline{n}^j_{i}-\underline{n}^j_i +1)$
(with convention the sum vanishes when $\tau_j=1$) is independent of
paths in $A_j({\bf U,V})$ being some combination of $\{U_i: 1\leq
i\leq j-1\}$.
Hence, with $k=j$ in
(\ref{start-end-bounds}), we observe $\alpha_j = \overline{n}_{\tau_j}^j+1$ takes on at most
$\lfloor 2n\delta +1\rfloor$ distinct values. The result now follows as paths in
$A_j({\bf U,V})$ are in $1:1$ correspondence with last switch times
$\alpha_j$. \qed

 {\bf Acknowledgement.}  We thank M. Bal\'azs and J. Pitman for
 helpful communications.

\bibliographystyle{plain}

\end{document}